\documentclass[10pt]{amsart}

\usepackage{latexsym, amscd, amsmath, verbatim, amssymb, graphics}

\newcommand{\numberseries}{\mdseries}   

\newlength{\thmtopspace}                
\newlength{\thmbotspace}                
\newlength{\thmheadspace}               
\newlength{\thmindent}                  

\newtheoremstyle{bfupright head,slanted body}
                {\thmtopspace}{\thmbotspace}
                {\slshape}{\thmindent}{\bfseries}{.}{\thmheadspace}
                {{\numberseries \thmnumber{(#2) }}\thmnote{#3}}

\newtheoremstyle{bfupright head,upright body}
                {\thmtopspace}{\thmbotspace}
                {\upshape}{\thmindent}{\bfseries}{.}{\thmheadspace}
                {{\numberseries \thmnumber{(#2) }}\thmnote{#3}}

\newtheoremstyle{bfit head,upright body}
                {\thmtopspace}{\thmbotspace}
                {\upshape}{\thmindent}{\upshape}{.}{\thmheadspace}
                {{\numberseries\thmnumber{(#2) }}
                {\bfseries\itshape\thmnote{\negthickspace#3}}}

\newtheoremstyle{it head,upright body}
                {\thmtopspace}{\thmbotspace}
                {\upshape}{\thmindent}{\upshape}{.}{\thmheadspace}
                {{\numberseries\thmnumber{(#2) }}
                {\itshape\thmnote{\negthickspace#3}}}

\newtheoremstyle{fixed bf head,slanted body}
                {\thmtopspace}{\thmbotspace}{\slshape}
                {\thmindent}{\bfseries}{.}{\thmheadspace}
                {{\numberseries \thmnumber{(#2) }}\thmname{#1}\thmnote{ (#3)}}

\newtheoremstyle{fixed bf head,upright body}
                {\thmtopspace}{\thmbotspace}{\upshape}
                {\thmindent}{\bfseries}{.}{\thmheadspace}
                {{\numberseries \thmnumber{(#2) }}\thmname{#1}\thmnote{ (#3)}}

\newtheoremstyle{independent paragraph}
                {\thmtopspace}{\thmbotspace}
                {\upshape}{\thmindent}{\upshape}{}{0pt}
                {\thmnote{#3 }}

\newtheoremstyle{subparagraph}
                {\thmbotspace}{\thmbotspace}
                {\upshape}{\thmindent}{\upshape}{}{0pt}
                {\thmnote{#3 }}

\newtheoremstyle{notes}
                {\thmtopspace}{\thmbotspace}
                {\ttfamily}{\thmindent}{\ttfamily\small }{}{0pt}
                {\thmnote{#3 }}

\theoremstyle{bfupright head,slanted body}
\newtheorem{res}{}[section]             \newtheorem*{res*}{}

\theoremstyle{bfit head,upright body}
                 \newtheorem*{com*}{}

\theoremstyle{bfupright head,upright body}
\newtheorem{bfhpg}[res]{}               \newtheorem*{bfhpg*}{}

\theoremstyle{it head,upright body}
               \newtheorem*{ithpg*}{}

\theoremstyle{fixed bf head,slanted body}
\newtheorem{thm}[res]{Theorem}          \newtheorem*{thm*}{Theorem}
\newtheorem{prp}[res]{Proposition}      \newtheorem*{prp*}{Proposition}
        \newtheorem*{cor*}{Corollary}
\newtheorem{lem}[res]{Lemma}            \newtheorem*{lem*}{Lemma}

\theoremstyle{fixed bf head,upright body}
       \newtheorem*{dfn*}{Definition}
      \newtheorem*{obs*}{Observation}
\newtheorem{rmk}[res]{Remark}           \newtheorem*{rmk*}{Remark}
\newtheorem{exa}[res]{Example}          \newtheorem*{exa*}{Example}
         \newtheorem*{exe*}{Exercise}
            \newtheorem{stp*}{Setup}
            \newtheorem*{que*}{Question}

\theoremstyle{independent paragraph}

\theoremstyle{subparagraph}

\theoremstyle{notes}

\newlength{\thmlistleft}        
\newlength{\thmlistright}       
\newlength{\thmlistpartopsep}   
\newlength{\thmlisttopsep}      
\newlength{\thmlistparsep}      
\newlength{\thmlistitemsep}     

\newcounter{eqc} 
  {\end{list}}%

\newcounter{prt}
  {\end{list}}%

\newcounter{rqm}
  {\end{list}}%

  {\end{list}}%

  {\end{list}}%

\newcommand{\pgref}[1]{(\ref{#1})}

\renewcommand{\eqref}[1]{\pgref{eq:#1}}

\renewcommand{\l}{\ell}
\newcommand{\ann}{\operatorname{ann}}
\newcommand{\depth}{\operatorname{depth}}
\newcommand{\h}{\operatorname{H}}
\newcommand{\rt}{\operatorname{reltype}}
\newcommand{\ar}{\operatorname{\mathsf{ar}}}

\begin{document}

\title{Strong Artin-Rees Property in Rings of dimension one and two}
\author{Janet Striuli}
\email{jstriuli2@math.unl.edu}
\address{Department of Mathematics, University of Nebraska-Lincoln, NE, 68588}
\maketitle
\begin{abstract}In this note we give a simple proof of the fact
that local  rings of dimension one have the  strong uniform Artin-Rees
property. Moreover, we give two examples of rings of dimension two
where the property fails.
\end{abstract}
\section{ Introduction}

In this paper   $(R,\mathbf{m} ,\mathsf{k})$  denotes a local Noetherian ring, and all  modules are finitely generated.

Let $I$ be an ideal  of $R$, let $M$ be an $R$-module and   $N$ a submodule of $M$. The Artin-Rees lemma states that there exists an integer $h$ depending on $I$, $M$ and $N$ such that  one has
\begin{eqnarray}\label{strong UAR}
I^nM \cap N=I^{n-h}(I^hM\cap N), \quad \text{ for all $n\geq h$}.
\end{eqnarray}
A weaker property, which is often the one used in the applications, is
\begin{eqnarray}\label{weak UAR}
I^{n}M\cap N\subset I^{n-h}N, \quad\text{ for all $n\geq h$} .
\end{eqnarray}
Much work has been done to determine whether $h$ can be chosen uniformly, in the sense that (\ref{strong UAR}) or (\ref{weak UAR})  would be satisfied simultaneously for every ideal belonging to a given  family, (\cite{Dun-Car-97}, \cite{Hune-92}, \cite{O'Car-87}, \cite{O'Car-91}, \cite{Plan-00}, \cite{Wang-97}). 
In  particular in \cite{Plan-00} it is proved that  there exists an integer $h$ such that (\ref{strong UAR}) is  satisfied simultaneously for every ideal if the ring is excellent of dimension one. In this note we want to give a simpler proof of this fact for local rings, with particular attention on the integer $h$. 

We also give an example of a family of ideals in a two dimensional ring for which there exists no integer $h$ such that (\ref{strong UAR}) holds for all ideals in the family.

\section {Strong Artin-Rees property in rings of dimension one  }
Let $R$ be a Noetherian ring  and let $\mathcal{W}$ be a family of ideals of $R$. Let $M$ be an $R$-module and let $N\subset M$ be a submodule. Let $h$ be an integer. Following the definitions in \cite{Hune-92}, we say that the pair $(N,M)$ has the \textit{strong Artin-Rees property with respect to $\mathcal{W}$ with Artin-Rees number $h$} if (\ref{strong UAR}) holds for all $I\in \mathcal{W}$. Notice that every  integer bigger than $h$ is an Artin-Rees number with respect to $\mathcal{W}$ for the pair $(N,M)$. We denote by $\ar_R(N,M;\mathcal{W})$ the least of such integers.

When  $\mathcal{W}$ is the family of all ideals, we say that  the pair $(N,M)$ has the \textit{strong Artin-Rees property} and  denote by $\ar_R(N,M)$ the  least of the Artin-Rees numbers.

We first show that it is enough to study the strong Artin-Rees property with respect to the family of $\mathbf{m}$-primary ideals. For this, we first need a lemma.

\begin{lem}\label{first}Let $M$ be an $R$-module. If  $N_{1},
N_{2}$ are two submodules of $M$ then there exists
$h=h((N_{1}+N_{2})\subseteq M)$ such that for every $n > h$ one
has 
\[N_{1} \cap (N_{2} +\mathbf{m}^{n}M)\subseteq  (N_{1} \cap N_{2}) + \mathbf{m}^{n-h}N_{1}.
\]
\end{lem}
\noindent\textit{Proof}.
By the Artin-Rees Lemma there exists $h$  such that for every $n > h$ we have $\mathbf{m}^nM\cap
(N_1+N_2)=\mathbf{m}^{n-h}(\mathbf{m}^{h}M \cap (N_1+N_2))\subset
\mathbf{m}^{n-h}(N_1+N_2)$. Then the following holds for $n>h$:
\begin{xxalignat}{3}
& &N_{1} \cap (N_{2} +\mathbf{m}^{n}M)&= N_{1}  \cap ( N_{2} + (\mathbf{m}^{n}M \cap
(N_{1}+N_{2})) \\
&\phantom{\square} & &\subseteq  N_{1} \cap(N_{2}+\mathbf{m}^{n-h}(\mathbf{m}^{h}M \cap (N_1+N_2)))\\
&\phantom{\square} & &\subseteq  N_{1} \cap (N_{2} + \mathbf{m}^{n-h}(N_{1}+N_{2})) \\
&\phantom{\square}& &\subseteq  N_{1} \cap (N_{2} + \mathbf{m}^{n-h}N_{1})\\
&\phantom{\square}& &\subseteq  N_{1} \cap N_{2} + \mathbf{m}^{n-h}N_{1} .& & \qed
\end{xxalignat}

\begin{rmk}\label{remarkfirst}Notice that if $h$ is an integer that satisfies Lemma
\ref{first}, then every  bigger integer does as well.
\end{rmk}
\begin{prp} Let $M$ be an $R$-module and $N\subset M$ a submodule. Let $\mathcal{W}$ be the family of $\mathbf{m}$-primary ideals.  Assume that $(N,M)$ has the strong uniform Artin-Rees property with respect to $\mathcal{W}$ then $\ar_R(N, M)\leq \ar_R(N,M;\mathcal{W})$.
\end{prp}
\begin{proof}
Let $h_0=\ar(N,M;\mathcal{W})$ and assume  by contradiction
that there exists $ I \subset R$ and  $ n \geq h_0 $ such that
$I^{n-h_0} (I^{h_0}M \cap N)\neq I^{n}M \cap N$.

On the other hand, for all $h>>0$ and for such a fixed $n$ and $h_0$, we have:
\begin{alignat*}{2}
 I^{n}M \cap N    & \subseteq (I +\mathbf{m}^{h})^{n}M \cap N, & \\
& \subseteq  (I+\mathbf{m}^{h})^{n-h_0}( (I +\mathbf{m}^{h})^{h_0}M \cap N),&\quad\text{by the definition of $h_0$}, \\
& \subseteq  I^{n-h_0}( (I+\mathbf{m}^{h})^{h_0}M \cap N) +\mathbf{m}^{h}M, &\quad \text{by expanding the powers},\\
& \subseteq  I^{n-h_0}((I^{h_0}+\mathbf{m}^{h})M \cap N) +\mathbf{m}^{h}M,&\\
& =  I^{n-h_0}((I^{h_0}M+\mathbf{m}^{h}M) \cap N) +\mathbf{m}^{h}M. &
\end{alignat*}
Let $h_1$ be an integer depending on $(I^{h_0}M+N) \subseteq M$
that satisfies Lemma \ref{first} with $N_{1}=N$, $N_{2}=I^{h_0}M$. By
Remark \ref{remarkfirst}, we may assume
$h_1 \geq n-h_0$. Applying Lemma \ref{first}, we have $(I^{h_0}M+\mathbf{m}^{h}M) \cap
N \subseteq  (I^{h_0}M \cap N) +\mathbf{m}^{h-h_1}M$, for every $h>h_1$. Therefore, the following holds:
\begin{align*}
 I^{n-h_0}((I^{h_0}M+\mathbf{m}^{h}M) \cap N) & +\mathbf{m}^{h}M \\
& \subseteq  I^{n-h_0}(I^{h_0}M \cap N + \mathbf{m}^{h-h_1}M) +\mathbf{m}^{h}M,\\
& \subseteq  I^{n-h_0}(I^{h_0}M \cap N ) + \mathbf{m}^{h-h_1+n-h_0}M +\mathbf{m}^hM,\\
&\subseteq I^{n-h_0}(  I^{h_0}M \cap N) + \mathbf{m}^{h-h_1+n-h_0}M.
\end{align*}
Putting together the right and the left end of the chain of inclusions, we  obtain that  $I^{n}M \cap N \subseteq I^{n-h_0}(  I^{h_0}M \cap N) +
\mathbf{m}^{h-h_1+n-h_0}M$, for every $h>h_1$. By taking the
intersection of the right side of the inclusion
over all $h>h_1$, we can conclude $I^{n}M \cap N \subseteq I^{n-h_0}(
I^{h_0}M \cap N)$. Since the reverse inclusion always holds, we conclude $I^{n-h_0} (I^{h_0}M \cap N)= I^{n}M \cap N$, contradicting the assumption.
\end{proof}

We  also need another kind of  reduction, see for example \cite[(2.4)]{Hune-92}.
\begin{lem}\label{faithfully}
Let $R \rightarrow S$ be a faithfully flat extension. Let $M$ be an $R$-module and $N\subset M$ a submodule. If $(N\otimes_R S,M\otimes_R S)$ has the strong uniform Artin-Rees property then  $\ar_R(N\subseteq M)\leq \ar_S(N\otimes_R S\subseteq M \otimes_R S)$.
\end{lem}

\begin{proof}A faithfully flat extension commutes with intersections.
\end{proof}

\begin{rmk}\label{infinite}Let $(R,\mathbf{m},\mathsf{k})$ be a local Noetherian
ring. Then $R \longrightarrow R[x]_{\mathbf{m}R[x]}$ is a faithfully flat extension
and $R[x]_{\mathbf{m}R[x]}$ has an infinite residue field.
\end{rmk}

\begin{prp}\label{second} 
Suppose $(R,\mathbf{m},\mathsf{k})$ is a  one-dimensional local Noetherian ring with infinite residue
field. Then there exists an integer $r=r(R)$,  such that for
every $\mathbf{m}$-primary ideal $I$ there exists $y \in I$ so that
$I^{n}=yI^{n-1}$, for every $n \geq h$. 
\end{prp}

\noindent \textit{Proof}.
First suppose that $R$ is Cohen-Macaulay and let $e$  be the multiplicity of
the ring. By \cite[Chapter 3,(1.1)]{Sally},  we
have that $\mu (I) \leq e$, where $\mu (I)$ denotes the minimal number
of generators of $I$ and $I$ is every $\mathbf{m}$-primary
ideal. Therefore,  $\mu (I^{e}) \leq e <
e+1$.  Hence, by \cite[Chapter 2, (2.3)]{Sally}, there exists $y \in I$
such that $I^{e}=yI^{e-1}$, so that for every $n \geq e$ we have
$I^{n}=yI^{n-1}$. Set $h$ to be $e$.

Next suppose $\depth (R)=0$, and let $0=q_{1} \cap q_{2}\dots \cap q_{s+1}$
be a minimal primary decomposition of $0$ where $q_{s+1}$ is  $\mathbf{m}$-primary
and set  $J=q_{1} \cap q_{2}\dots \cap q_{s}$. Then $R/J$ is
Cohen-Macaulay and there exists a $h_{0}$ such that $\mathbf{m}^{h_{0}}J=0$. Let $e_{1} $ be the multiplicity of $R/J$ then, by the
above case, there exists a $y \in I$ such that for every $n \geq e_{1}$
we have $I^{n} \subseteq yI^{n-1}+J$ and hence $I^{n}\subseteq yI^{n-1}+I^n\cap J$, for every $n>e_1$. By \cite[(4.2)]{Hune-92},
there exists a $ h_{1}$, depending just on $R$ and $J$ such that for
every $n \geq h_{1}$ and  every ideal $I \subset R$ we have $I^{n} \cap
J \subset I^{n-h_1}J$. Hence, for every $n\geq h= \max \{e_{1},h_{0}+h_1\}$ one has the following inclusions:
\begin{xxalignat}{3}
& &I^{n} &\subseteq yI^{n-1}+I^{n} \subseteq yI^{n-1}+I^{n}\cap J \\
&\phantom{square}& &\subseteq yI^{n-1}+I^{n-h_1}J\subseteq yI^{n-1}+\mathbf{m}^{h_0}J = yI^{n-1}. & &\qed
\end{xxalignat}
We are now ready to prove the main theorem
If $M$ is a finite length module we denote by $\l(M)$ its length.
\begin{thm}\label{main} Let $(R,m,\mathsf{k})$ be a one-dimensional local
ring with infinite residue field. Then $R$ has the strong uniform Artin-Rees property.

Moreover, if $M$ is  a finitely generated $R$-module and  $N\subset M$ a submodule, then 
\[
\ar_R(N, M)\leq \max\{r,\l(\h^0_\mathbf{m}(M/N))\}+\l(\h^0_\mathbf{m}(M/N)),
\]
 where  $r=r(R)$ is an integer as in Proposition \ref{second}.
\end{thm}
\begin{proof}  Let $I$ be an $\mathbf{m}$-primary ideal. Set 
\begin{equation*}
h_1 = \l(\h^0_\mathbf{m}(M/N)) \quad \text{and}\quad h= \max\{r,\l(\h^0_\mathbf{m}(M/N))\}+\l(\h^0_\mathbf{m}(M/N)).
\end{equation*}
Assume first that $M/N$ is Cohen-Macaulay. By Proposition
\ref{second} we can choose $y \in I$ such that $y$ is a
non-zerodivisor in $M/N$, so that for $n > h=r$,
\begin{eqnarray*}
I^{n}M \cap N &=& yI^{n-1}M \cap N, \\
     & \subseteq & y(I^{n-1}M\cap N), \quad \text{by the property of $y$},\\
& \subseteq & I(I^{n-1}M \cap N), \quad \text{since $y \in I$}.
\end{eqnarray*}

Now suppose that $M/N$ is not Cohen-Macaulay and let $M^{'}/N =H_{m}^{0}(M/N)$.  For every $n \geq h$ and every $I \subseteq R$ we have:
\begin{alignat*}{2}
I^{n}M \cap N &= I^{n}M \cap M^{'} \cap N,& \quad \text{since $N
\subset M'$},\\
&= I^{n-r}(I^{r}M \cap M^{'}) \cap N,& \quad \text{since $M/M'$ is Cohen-Macaulay},\\
&\subseteq I^{n-r}(I^{r}M \cap M^{'}),&\quad \text{since $ n-r>h_1$ and $I^{h_1}M'\subset N$},\\
&= I^{n-r-h_1}I^{h_1}(I^{r}M \cap M^{'}),&\\
&=I^{n -h}(I^{h_1}(I^{r}M\cap M')\cap N),&\quad \text{since $ I^{h_1}M'
\subset N$},\\
&\subseteq I^{n-h}(I^{r+h_1}M\cap M' \cap N),&\quad\\
&=I^{n-h}(I^{h}M\cap N).&
\end{alignat*}
proving the theorem.
\end{proof}

\begin{bfhpg}[Relation Type]\label{Planas}
Let $I=(f_{1}, \dots, f_{n})$ be an ideal in $R$. Map the  polynomial
ring, with the standard grading, $R[x_{1}, \dots ,x_{n}]$ onto the Rees algebra $R[It]$ by sending
$f_{i}$ to $x_{i}t$. Let $L$ be the kernel of this map. Then $L$ is an
homogeneous ideal and
the relation type of $I$ is defined to be the minimum integer $h$ such that
the ideal $L$ can be generated by elements of degree less or equal
than $h$. It is denoted by $\rt(I)$.  This number does not depend on the choice of the minimal generators
of the ideal $I$.

 If $(R,\mathbf{m},\mathsf{k})$ is a one-dimensional, Cohen-Macaulay local
ring and $I$ is an \linebreak $\mathbf{m}$-primary ideal, then $\rt(I) \leq e$, where $e$
is the multiplicity of $R$, see  \cite{AchSch-82}.
\end{bfhpg}
The following lemma had been proved by Wang in \cite{Wang-97B} for parameters ideals. The same argument applies for every ideal, we include it here for simplicity.
\begin{lem}\label{relationtype} Let $(R,\mathbf{m},\mathsf{k})$ be a local ring and $J$ be an ideal of
$R$; denote $\bar{R}=R/J$. Let $I=(x_{1}, \dots ,x_{m})$ be an ideal
of $R$ and suppose that $\rt (I \bar{R}) \leq h$, for some $h>0$. Then
for every $n > h$,
\[
I^n \cap J =I^{n-h}(I^{h} \cap J). 
\]
\end{lem}
\begin{proof}Let $n> h$ and let $x \in I^{n} \cap J$. Then there
exists a polynomial $F$ in $R[X_{1}, \dots
,X_{m}]$, homogeneous of degree $n$, such that $F(x_{1}, \dots
,x_{m})=x$. Modulo $J$, $\bar{F}$ is a relation on the
$\bar{x_{i}}$'s, so by hypothesis there are polynomials $G_{i}$ of
degree $h$, and
$H_{i}$, of degree $n-h$, such that $\bar{F}=\sum \bar{G_{i}} \bar{H_{i}}$ in $
\bar{R}[X_{1}, \dots ,X_{m}]$ and $\bar{G_{i}}$ are relations on the
$\bar{x_{i}}$. Therefore, $F= \sum G_{i}H_{i} +K$ for some $K \in
R[X_{1},\dots ,X_{m}]$ of degree $n$ and coefficients in $J$. Since:
\begin{eqnarray*}
K(x_{1}, \dots ,x_{m}) &\in& JI^n \subset I^{n-h}(I^{h}\cap J),\\
G_{i}(x_{1}, \dots ,x_{m}) &\in& I^h \cap J, \quad \text{and}\\
H_{i}(x_{1} , \dots ,x_{m}) &\in& I^{n-h},
\end{eqnarray*}
$x=F(x_{1}, \dots,x_{m})=\sum G_{i}(x_{1},\dots ,x_{m})H_{i}(x_{1},
\dots ,x_{m}) \in I^{n-h}(I^{h} \cap J)$.
\end{proof}

\begin{lem} \label{same}Let $(R,\mathbf{m},\mathsf{k})$ a Noetherian local ring. If $J$ is an
ideal of $R$ such that $\dim(R/J)\leq 1$ then $(J,R)$ has the strong Artin-Rees property.
\end{lem}
\noindent\textit{Proof}. If $\dim(R/J)=0$ then there exists a power of the maximal
ideal $\mathbf{m}^h\subset J$. Therefore, for $n> h$ and for every ideal
$I$ we have the following:
\begin{eqnarray*}
I^{n}\cap J = I^n=II^{n-1}=I(I^{n-1}\cap J).
\end{eqnarray*}
Assume $\dim(R/J)=1$. By Lemma \ref{first} it is enough to show that  $(J,R)$ has the strong Artin-Rees property with respect to the family of $\mathbf{m}$-primary ideals.
Suppose that $R/J$ is Cohen-Macaulay, then the conclusion holds by \ref{Planas}  and by
Lemma \ref{relationtype}.\\
Suppose $R/J$ has dimension one and it is not Cohen-Macaulay. Let $J
\subset J'$ such that $R/J'$ is Cohen-Macaulay and let $h_0$ such that
$\mathbf{m}^{h_0}J'\subset J$. By the Cohen-Macaulay case there exists
an Artin-Rees number $h_1=h_1(J' \subset R)$. We may assume $h_1>h_0$. Let 
$h=h_1+h_0 $.  Then, with an argument we already used, for every $n>h$
one has 
\begin{xxalignat}{3}
& & I^{n} \cap J &= I^{n} \cap J' \cap J \\
&\phantom{square}&            &= I^{n-h_1}( I^{h_1}\cap J')\cap J\\
&\phantom{square}&            &= I^{n-h_1}( I^{h_1}\cap J')\\
&\phantom{square}&             &= I^{n-h_1-h_0}I^{h_0}( I^{h_1}\cap J')\\
&\phantom{square}&             &= I^{n-h_1-h_0}(I^{h_0}( I^{h_1}\cap J')\cap J)\\
&\phantom{square}&             &\subseteq I^{n-h}(I^{h} \cap J' \cap J)\\
&\phantom{square}&             &= I^{n-h}(I^{h} \cap J). & &\qed
\end{xxalignat}

\begin{prp}\label{module->ring}Let $(R,\mathbf{m},\mathsf{k})$ be a local Noetherian ring. Let $M$ be an \mbox{$R$-module} and $N\subseteq M$ a submodule. Let $J\subset\ann(M/N)$ be an ideal of $R$. If  $(J, R)$  and $(N/JM,M/JM)$ have the strong uniform Artin-Rees property, then 
\[
\ar_R(N,M)\leq \max\{\ar_R(J,R),\ar_{R/J}(N/JM,M/JM)\}.
\]
In particular if $\dim(M/N)=1$ and  the residue field is infinite then 
\[
\ar_R(N,M)\leq \max\{\ar_R(J,R),\max\{r(R/J), \l(\h^0_\mathbf{m}(M/N)\}+\l(\h^0_\mathbf{m}(M/N)\}.
\]
\end{prp}
\noindent \textit{Proof.}
The second statement follows from the first and Theorem \ref{main}. 
For the first part, let $h=\max\{\ar_R(J,R),\ar_{R/J}(N/JM,M/JM)\}$.
Let $\phi: R^{m} \rightarrow M$, a surjection of a free module onto
$M$. Denote by $K=\ker (\phi)$ and by $L=\phi ^{-1}(N)$, the pre-image
of the submodule $N\subset M$.  Then, as shown in
\cite{Dun-Car-97}, it is enough to show that there exists a $h$
such that for every $n>h$ and for every ideal $I\subset R$, we have
$I^{n}R^{m} \cap L = I^{n-h}(I^{h}R^{m} \cap L)$. Therefore, without loss of
generality we may assume $M$ is a free module. 

Since $h\geq\ar_{R/J}(N/JM,M/JM)$, for every
$n>h $ and for every ideal $I$, we have $I^{n}M\cap N \subset
I^{n-h}(I^{h}M \cap N) + JM$. Therefore, 
\[
I^{n}M\cap N \subset
I^{n-h}(I^{h}M \cap N) + JM \cap I^{n}M = I^{n-h}(I^{h}M \cap N) +
(I^{n}\cap J)M,
\]
 where the last equality holds since $M$ is a free
module. Since $h\geq\ar_R(J,R)$, we have $I^{n}\cap J = I^{n-h}(I^{h}\cap J)$. Hence,
\begin{xxalignat}{3}
& & I^{n}M\cap N &=I^{n-h}(I^{h}M \cap N)+I^{n-h}(I^{h}\cap
J)M\\
&\phantom{square}& &=I^{n-h}(I^{h}M \cap N)+I^{n-h}(I^{h}M\cap JM)\\
&\phantom{square}& &\subset
I^{n-h}(I^{h}M \cap N), \quad \text{since $JM \subseteq N$}.& & \qed
\end{xxalignat}

\section{Dimension two}
The following example (see \cite{Wang-97}), shows that the uniform Artin-Rees property does
not hold for two dimensional rings.
\begin{exa}Let $R=\mathsf{k}[x,y,z]/(z^{2})$. Consider the following family of
ideals:
\[
I_n=(x^n,y^n, x^{n-1}y +z),
\]
for every $n \in \mathbb{N}$. Let $J$ the ideal generated by $z$. 

We want to show that 
$I_{n}(I^{n-1}_{n} \cap J)\neq I^{n}_{n}\cap J$, for every $n\geq 2$. In
particular we will show that 
\[
 x^{(n-1)^2}y^{n-1}z \in  I^{n}_{n}\cap J \quad \text{but} \quad x^{(n-1)^2}y^{n-1}z \notin I_{n}(I^{n-1}_{n} \cap J).
\]
Denote  $ x^{(n-1)^2}y^{n-1}z$ by $\xi$. 

The ideal $I_n$
is a homogeneous ideal if we assign degree  one to $x$ and $y$ and
degree $n$ to $z$. With such grading $\xi$ has degree $(n-1)^2+n-1 +n= n^2$.  Since
$x^{(n-1)^2}y^{n-1}z = (x^{n-1}y+z)^{n}-(x^{n})^{n-1}y^{n} \in
I_{n}^{n}$ the first claim holds. 

Suppose  $x^{(n-1)^2}y^{n-1}z
\in I_{n}( I^{n-1}_{n}\cap J)$, this remains true  modulo   $(x^{(n-1)^2+1}, y^n)$. The ideal $I^{n-1}_{n}$ modulo $ (x^{(n-1)^2+1},y^n)R$ is generated
by 
\[
\{x^{n(n-1-i)}(x^{n-1}y+z)^{i} \mid i=0,1,\dots, n-1\}.
\]
Moreover,
\begin{eqnarray*}
x^{n(n-1-i)}(x^{n-1}y+z)^{i} &=& x^{n(n-1-i)}(x^{(n-1)i}y^i
+x^{(n-1)(i-1)}y^{i-1}z)\\
&=& x^{n^2-n-i}y^{i}+x^{n^2-2n-i+1}y^{i-1}z.
\end{eqnarray*}
But $n^2-n-i\geq (n-1)^2+1$ for $i\leq n-2$. Therefore, $I_{n-1}^{n}$ modulo
$(x^{(n-1)^2+1},y^n)$ is generated by
\[
\{x^{(n-1)^2}y^{n-1}+x^{(n-1)(n-2)}y^{(n-2)}z ,\quad
x^{n^2-2n-i+1}y^{i-1}z \mid i=1,\dots
,n-2\}.
\]
Let 
\begin{eqnarray*}
f= x^{(n-1)^2}y^{n-1}+x^{(n-1)(n-2)}y^{(n-2)}z,\\
g_{i}=x^{n^2-2n-i+1}y^{i-1}z.
\end{eqnarray*}

Let  $hf+\sum
h_ig_i$ be a homogeneous element of $I^{n-1}_{n}\cap J$ that appear in the expression on
$\xi$ as element of  $I_n(I^{n-1}_{n}\cap J)$. By degree reasons we can assume $h$ is not a constant polynomial. 

 Let $m(x,y,z)$ be a homogeneous
monomial of $h$. If $z$ does not divide $m$, then
\begin{align*}
& m(x,y,z)f=m'(x,y)x^{(n-1)(n-2)+1}y^{(n-2)}z\\
 \text{or} \quad & m(x,y,z)f=m'(x,y)x^{(n-1)(n-2)}y^{(n-2)+1}z;
\end{align*}
 if $z$ does divide $m$ then
$m(x,y,z)f=m'(x,y)x^{(n-1)^2}y^{n-1}z$, with $m'$ possibly a unit.
By a degree counting we can see that $\deg(hf)\geq
n^2-n+1$. Therefore, for every element $a\in I_{n-1}$ we have
$\deg (ahf)>n^2 =\deg(\xi)$. This shows a contradiction.
\end{exa}

The following example shows that the Artin-Rees property fails in a
two dimensional ring, even if the ring is reduced.
\begin{exa}Let $R=\mathsf{k}[x,y,z]/xz$. Consider the following family of
ideals:
\begin{displaymath}
I_{n}=(x^n,y^n,x^{n-1}y+z^n),
\end{displaymath}
for every $n\in \mathbb{N}$. Let $J=(z)$. Again, we claim that 
$I_{n}(I^{n-1}_{n} \cap J)\neq I^{n}_{n}\cap J$ for every $n\geq 1$. We will
show that 
\begin{displaymath}
z^{n^2} \in  I^{n}_{n}\cap J \quad \text{ but} \quad  z^{n^2} \notin
I_{n}(I^{n-1}_{n} \cap J).
\end{displaymath}
 Indeed,  $z^{n^2}=(x^{n-1}y+z^n)^{n}-
(x^{n})^{n-1}y^{n} \in I_{n}^{n}$ and trivially $z^{n^2}\in J$.

On the other hand $I^{n-1}_{n}$ is generated by:
\begin{eqnarray*}
 \{x^{n(n-1)}, x^{(n-1)^2}y^{n-1}+z^{n(n-1)}, y^{n}L , x^{(n-1)^2+i}y^{n-1-i} \mid i =1,\dots n-1\},
\end{eqnarray*}
for some ideal $L$ in $R$. Notice that if $z^{n^2} \in
I_n(I_{n}^{n-1}\cap J)$ then this also holds  modulo $y^n$. Moreover, if a homogeneous
element 
\[f(x,y)x^{n(n-1)}+g(x,y,z)(x^{(n-1)^2}y^{n-1}+z^{n(n-1)})+\sum_{i=1}^{n-1} h_{i}(x,y)
x^{(n-1)^2+i}y^{n-1-i}
\]
 is in $J$, writing $g(x,y,z)=g''(x,y)+zg'(x,y,z)$, we
see that
\begin{displaymath} 
f(x,y)x^{n(n-1)}+g''(x,y)x^{(n-1)^2}y^{n-1}+\sum h_{i}(x,y)
x^{(n-1)^2+i}y^{n-1-i}=0.
\end{displaymath}
But if this is the case, since $xz=0$ in $R$, we have 
\begin{eqnarray*}
fx^{n(n-1)}+g(x^{(n-1)^2}y^{n-1}+z^{n(n-1)})+\sum h_{i}
x^{(n-1)^2+i}=zg'z^{n(n-1)}.
\end{eqnarray*}
 But $zg'z^{n(n-1)}$
is an homogeneous element of degree at least $n^2 -n +1$ and 
multiplication by any element in $I_{n}$ increases the degree by
$n$. Therefore, any element in $I_{n}(I_{n}^{n}\cap J)$ has degree at
least $n^2+1$ while $z^{n^2}$ has degree strictly smaller.
\end{exa}

\bibliographystyle{amsplain}

\begin{thebibliography}{1}

\bibitem{AchSch-82}
R{\"u}diger Achilles and Peter Schenzel, \emph{A degree bound for the defining
  equations of one-dimensional tangent cones}, Seminar D. Eisenbud/B. Singh/W.
  Vogel, Vol. 2, Teubner-Texte zur Math., vol.~48, Teubner, Leipzig, 1982,
  pp.~77--87. \MR{686461 (84g:14027)}

\bibitem{Dun-Car-97}
A.~J. Duncan and L.~O'Carroll, \emph{A full uniform {A}rtin-{R}ees theorem}, J.
  Reine Angew. Math. \textbf{394} (1989), 203--207. \MR{977443 (90c:13011)}

\bibitem{Hune-92}
Craig Huneke, \emph{Uniform bounds in {N}oetherian rings}, Invent. Math.
  \textbf{107} (1992), no.~1, 203--223. \MR{93b:13027}

\bibitem{O'Car-91}
L.~O'Carroll, \emph{A note on {A}rtin-{R}ees numbers}, Bull. London Math. Soc.
  \textbf{23} (1991), no.~3, 209--212. \MR{92i:13001a}

\bibitem{O'Car-87}
Liam O'Carroll, \emph{A uniform {A}rtin-{R}ees theorem and {Z}ariski's main
  lemma on holomorphic functions}, Invent. Math. \textbf{90} (1987), no.~3,
  647--652. \MR{914854 (89h:13013)}

\bibitem{Plan-00}
Francesc Planas-Vilanova, \emph{The strong uniform {A}rtin-{R}ees property in
  codimension one}, J. Reine Angew. Math. \textbf{527} (2000), 185--201.
  \MR{2001g:13051}

\bibitem{Sally}
Judith~D. Sally, \emph{Numbers of generators of ideals in local rings}, Marcel
  Dekker Inc., New York, 1978. \MR{0485852 (58 \#5654)}

\bibitem{Wang-97B}
Hsin-Ju Wang, \emph{The relation-type conjecture holds for rings with finite
  local cohomology}, Comm. Algebra \textbf{25} (1997), no.~3, 785--801.
  \MR{1433435 (98d:13017)}

\bibitem{Wang-97}
\bysame, \emph{Some uniform properties of {$2$}-dimensional local rings}, J.
  Algebra \textbf{188} (1997), no.~1, 1--15. \MR{98e:13006}

\end{thebibliography}
\providecommand{\bysame}{\leavevmode\hbox to3em{\hrulefill}\thinspace}
\providecommand{\MR}{\relax\ifhmode\unskip\space\fi MR }
\providecommand{\MRhref}[2]{%
  \href{http://www.ams.org/mathscinet-getitem?mr=#1}{#2}
}
\providecommand{\href}[2]{#2}

\end{document}